\documentclass[12pt]{amsart}

\usepackage{amsmath} \usepackage{amsthm} \usepackage{mathrsfs}
\usepackage{latexsym}
\usepackage{amssymb}
\usepackage{enumerate}

\newtheorem{thm}{Theorem}[section]
\newtheorem{lem}[thm]{Lemma}

\newtheorem{prop}[thm]{Proposition}

\newtheorem*{thm*}{Theorem}
\newtheorem*{cor*}{Corollary}
\newtheorem*{question*}{Question}

\theoremstyle{definition}
\newtheorem{defn}[thm]{Definition}

\newtheorem{example}[thm]{Example}

\theoremstyle{remark}

\newcommand\cfrak{\mathfrak{c}}
\newcommand\dfrak{\mathfrak{d}}

\newcommand\Dcal{\mathscr{D}}

\newcommand\Pcal{\mathscr{P}}

\newcommand\Rbb{\mathbb{R}}

\newcommand{\dom}{\operatorname{dom}}

\newcommand{\otp}{\operatorname{otp}}

\newcommand\triord{\triangleleft}

\newcommand\axiom{\mathrm}
\newcommand\MA{\axiom{MA}}

\newcommand\PFA{\axiom{PFA}}
\newcommand\BPFA{\axiom{BPFA}}

\newcommand\MM{\axiom{MM}}
\newcommand\BMM{\axiom{BMM}}
\newcommand\MRP{\axiom{MRP}}
\newcommand\SRP{\axiom{SRP}}

\newcommand\AC{\axiom{AC}}

\newcommand\hull{\operatorname{hull}}

\newcommand\fin{\mathrm{fin}}

\newcommand\NS{\mathrm{NS}}

\renewcommand{\>}{\rangle}

\renewcommand{\epsilon}{\varepsilon}

\newcommand\mand{\textrm{ and }}

\subjclass[2000]{03E05, 03E10, 03E47, 03E65}


\thanks{This paper served as a cornerstone for the proposal
for my NSF grant DMS--0401893;
revisions were made to the paper while I was supported by this grant.}

\title[Set mapping reflection]
{Set mapping reflection}

\author{Justin Tatch Moore}

\begin{document}

\begin{abstract}
In this note we will discuss a new reflection principle
which follows from the Proper Forcing Axiom.
The immediate purpose will be to prove that
the bounded form of the Proper Forcing Axiom implies both
that $2^{\omega} = \omega_2$ and that
$L(\Pcal(\omega_1))$ satisfies the Axiom of Choice.
It will also be demonstrated that this reflection principle
implies that $\square(\kappa)$ fails for all
regular $\kappa > \omega_1$.
\end{abstract}

\maketitle
\section{Introduction}
The notion of properness was introduced by Shelah and is a weakening of both
the countable chain condition and the property of being countably closed.
Its purpose was to provide a property of forcing notions which implies that they preserve
$\omega_1$ and which is preserved under countable support iterations.
With the help of a supercompact cardinal,
one can prove the consistency of the following statement (see \cite{Yorkshireman}).
\begin{description}

\item[$\PFA$]
If $\Pcal$ is a proper forcing notion and $\Dcal$ is a family of dense subsets of $\Pcal$
of size $\omega_1$ then there is a filter $G \subseteq \Pcal$ which meets every element
of $\Dcal$.

\end{description}
The Proper Forcing Axiom (PFA for short) is therefore a strengthening of the
better known and less technical $\MA_{\omega_1}$\cite{iterated_cohen_ext}.
It has been extremely useful, together with the stronger Martin's Maximum (MM) \cite{FMS},
in resolving questions left unresolved by Martin's Axiom.

Early on it was known that it was not possible to replace
$\omega_1$ by $\omega_2$
and get a consistent
statement (see \cite{PFA:Baumgartner}).
The stronger forcing axiom MM was known to already
imply that the continuum is $\omega_2$ \cite{FMS}.
Later Todor\v{c}evi\'{c} and Veli\v{c}kovi\'{c}
showed that $\PFA$ also implies that the continuum is $\omega_2$
(see \cite{comparing_continuum} and \cite{FA_stationary}).
This proof used a deep analysis of the gap structure of $\omega^\omega/\fin$
and of the behavior of the oscillation map.
Their proof, however, was less generous than some of the proofs that
the continuum was $\omega_2$ from Martin's Maximum.
In particular, while $\MM$ implies that $L(\Pcal(\omega_1))$ satisfies $\AC$
\cite{Pmax}, the same was not known for $\PFA$ (compare to the final remark
section 3 of \cite{FA_stationary}).

Martin's Maximum was also shown to have a variety of large cardinal consequences.
Many of these are laid out in \cite{FMS}.
Much of this was proved via stationary reflection principles which seemed to typify the
consequences of $\MM$ which do not follow from $\PFA$.
Todor\v{c}evi\'{c} showed that $\PFA$ implies that the combinatorial principle
$\square(\kappa)$ fails for all $\kappa > \omega_1$ (see \cite{set_theory:Bekkali}).
This, combined with modern techniques in inner model theory \cite{square_core},
gives a considerable lower bound on the consistency strength of $\PFA$.\footnote{
It is not known if ``$\square(\kappa)$ fails for all regular $\kappa > \omega_1$'' is
equiconsistent with the existence of a supercompact cardinal.}
Both the impact on the continuum and the large cardinal strength of these
forcing axioms have figured prominently in their development.

The purpose of this note is to introduce a new reflection principle,
$\MRP$, which follows from the Proper Forcing Axiom.
The reasons are threefold.
First, this axiom arose as a somewhat natural abstraction of one its
consequences which in turn implies that there is a
well ordering of $\Rbb$ which is $\Sigma_1$-definable over
$(H(\omega_2),\in)$.
A corollary of the proof will be that the Bounded Proper Forcing Axiom
implies that there is such a well ordering of $\Rbb$,
thus answering a question from the folklore
(see Question 35 of \cite{local_reflection}).

Second, this principle seems quite relevant in studying
consequences of the Proper Forcing Axiom which do not follow from
the $\omega$-Proper Forcing Axiom.
The notion of $\omega$-properness was introduced by Shelah in the
course of studying preservation theorems for not adding reals
in countable support iterations (see \cite{proper_forcing}).
For our purpose it is sufficient to know that both c.c.c. and countably
closed forcings are $\omega$-proper and that $\omega$-proper forcings are
preserved under countable support iterations.
While I am not aware of the $\omega$-$\PFA$ having been studied in
the literature,
nearly all of the studied consequences of $\PFA$ are actually consequences
of the weaker $\omega$-$\PFA$.\footnote{For example: $\MA_{\omega_1}$,
the non-existence of S-spaces \cite{forcing_partition},
all $\omega_1$-dense sets of reals are isomorphic \cite{reals_isomorphic},
the Open Coloring Axiom \cite{partition_problems},
the failure of $\square_\kappa$ for all regular $\kappa > \omega_1$ \cite{set_theory:Bekkali},
the non-existence of Kurepa trees \cite{PFA:Baumgartner}.}
It is my hope and optimism that the Mapping Reflection Principle will
be useful tool in studying the consequences of $\PFA$ which do not
follow from the $\omega$-$\PFA$ in much the same way that
the Strong Reflection Principle has succeeded in implying the typical consequences
of Martin's Maximum which do not follow from the Proper Forcing Axiom.

Finally, like the Open Coloring Axiom, the Ramsey theoretic formulation of
Martin's Axiom, and the Strong Reflection Principle,
this principle can be taken as a black box and used without knowledge of forcing.
The arguments using it tend to be rather elementary in nature and
require only some knowledge of the combinatorics of the club filter
on $[X]^\omega$ and L\"owenheim-Skolem arguments.

The main results of this note are summarized as follows.

\begin{thm}
The Proper Forcing Axiom implies the Mapping Reflection Principle.
\end{thm}

\begin{thm}
The Mapping Reflection Principle
implies that $2^{\omega} = 2^{\omega_1} = \omega_2$
and that $L(\Pcal(\omega_1))$ satisfies the Axiom of Choice.
\end{thm}

\begin{thm} \label{BPFA_c}
The Bounded Proper Forcing Axiom implies that
$2^{\omega} = \omega_2$ and that $L(\Pcal(\omega_1))$ satisfies
the Axiom of Choice.
\end{thm}

\begin{thm}
The Mapping Reflection Principle
implies that $\square(\kappa)$ fails for every regular
$\kappa > \omega_1$.
\end{thm}

The notation used in this paper is more or less standard.
If $\theta$ is a regular cardinal then $H(\theta)$ is the collection of all sets of
hereditary cardinality less than $\theta$.
As is common, when I refer to $H(\theta)$ as a structure I
will actually mean $(H(\theta),\in,\triord)$
where $\triord$ is some well order of $H(\theta)$ which can
be used to compute Skolem functions and hence generate the club
$E \subseteq [H(\theta)]^{\omega}$ of
countable elementary submodels of $H(\theta)$.
If $X$ is a set of ordinals then $\otp(X)$
represents the ordertype of $(X,\in)$
and $\pi_X$ is the unique collapsing isomorphism from $X$ to $\otp(X)$.
While an attempt has been made to keep parts of this paper self contained,
a knowledge of proper forcing is assumed in Section \ref{PFA:sec}.
The reader is referred to \cite{PFA:Baumgartner}, \cite{proper_forcing}, and
\cite{partition_problems} for more reading on proper forcing and $\PFA$.
Throughout the paper the reader is assumed to have a
familiarity with set theory
(\cite{set_theory:Jech} and \cite{set_theory:Kunen} are standard references). 

\section{The Mapping Reflection Principle}

The following definition will be central to our discussion.
Recall that for an uncountable set $X$,
$[X]^{\omega}$ is the collection of all countable subsets of $X$.
\begin{defn}
Let $X$ be an uncountable set, $M$ be a countable
elementary submodel of $H(\theta)$ for some
regular $\theta$ such that $[X]^{\omega} \in M$.
A subset $\Sigma$ of $[X]^{\omega}$ is \emph{$M$-stationary}
if whenever $E \subseteq [X]^{\omega}$ is a club in $M$
there is an $N$ in $E \cap \Sigma \cap M$.
\end{defn}
\begin{example}
If $M$ is a countable elementary submodel of $H(\omega_2)$ and
$A \subseteq M \cap \omega_1$ has order type less than $\delta = M\cap \omega_1$
then $\delta \setminus A$ is $M$-stationary.
\end{example}

The set $[X]^{\omega}$ is equipped with the Ellentuck
topology obtained by declaring the sets
$$[x,N] = \{Y \in [X]^\omega: x \subseteq Y \subseteq N\}$$
to be open for all $N$ in $[X]^{\omega}$ and finite $x \subseteq N$.
In this paper ``open'' will always refer to this topology.
It should be noted that the sets which are closed in the
Ellentuck topology and cofinal
in the order structure generate the closed unbounded filter on $[X]^{\omega}$.

For ease of reading I will make the following definition.
\begin{defn}
A set mapping $\Sigma$ is said to be \emph{open stationary} if, for some
uncountable set $X$ and regular cardinal $\theta$ with $X$ in $H(\theta)$,
it is the case that elements of the domain of $\Sigma$ are elementary submodels
of $H(\theta)$ which contain $M$ and $\Sigma(M) \subseteq [X]^\omega$
is open and $M$-stationary for all $M$ in the domain of $\Sigma$.
If necessary the underlying objects $X$ and $\theta$ will be referred
to as $X_\Sigma$ and $\theta_\Sigma$.
\end{defn}

The following is among the simplest example of an open stationary set mapping.
\begin{example}
Let $r:\omega_1 \to \omega_1$ be regressive on the limit ordinals.
If $\Sigma$ is defined by putting $\Sigma(P) = (r(\delta),\delta)$
for $P$ a countable elementary submodel of $H(\omega_2)$ then
$\Sigma$ is open and stationary.
\end{example}
This motivates the following reflection principle which asserts that
this example is present inside any open stationary set mapping.
\begin{description}

\item[$\MRP$]
If $\Sigma$ is an open stationary set mapping whose domain is a club then
there is a continuous $\in$-chain $\<N_\nu : \nu < \omega_1\>$ in
the domain of $\Sigma$ such that for all limit $0 < \nu < \omega_1$ there is a
$\nu_0 < \nu$ such that
$N_\xi \cap X_\Sigma \in \Sigma(N_\nu)$ whenever $\xi$ is
in the interval $(\nu_0,\nu)$.

\end{description}
We now continue with the first example.
\begin{example}
An immediate consequence of $\MRP$ is that if $C_\delta$ is a cofinal $\omega$-sequence
in $\delta$ for each countable limit ordinal $\delta$ then there is
a club $E \subseteq \omega_1$ such that $E \cap C_\delta$ is finite for all
$\delta$.\footnote{It is easy to verify, however, that this consequence
of $\MRP$ can not be forced with an $\omega$-proper forcing.
That this statement follows from $\PFA$ appears in \cite{proper_forcing}.
It is the only example in the literature that I am aware of which
is a combinatorial consequence of $\PFA$ but not of $\omega$-$\PFA$.}
\end{example}
We will see in the discussion below that, unlike in $[X]^{\omega}$,
there are non-trivial partitions of $[X]^\omega$ into open $M$-stationary sets if
$X$ has size at least $\omega_2$.

\section{$\PFA$ implies $\MRP$}
\label{PFA:sec}
The purpose of this section is to prove the following theorem.
Recall that a forcing notion $\Pcal$ is \emph{proper} if whenever
$M$ is a countable elementary submodel of $H(|2^{\Pcal}|^+)$ containing $\Pcal$
and $p$ is in $\Pcal \cap M$, there is a $\bar p \leq p$ which is \emph{$(M,P)$-generic}.\footnote{
Here condition $\bar p$ is $(M,\Pcal)$-generic if whenever $r$ is an extension of $\bar p$
and $D \subseteq \Pcal$ is a predense set in $M$, there is a $s$ in $D \cap M$ such that
$s$ is compatible with $r$.}

\begin{thm}
$\PFA$ implies $\MRP$.
\end{thm}

\begin{proof}
Let $\Sigma$ be a given open stationary set mapping defined on a club of models and
abbreviate $X = X_\Sigma$ and $\theta = \theta_\Sigma$.
Let $\Pcal_\Sigma$ denote the collection of all continuous $\in$-increasing maps
$p:\alpha + 1 \to \dom(\Sigma)$ where $\alpha$ is a countable ordinal such that
for all $0 < \nu \leq \alpha$ there is a $\nu_0 < \nu$ with
$p(\xi) \cap X \in \Sigma(p(\nu))$ whenever $\nu_0 < \xi < \nu$.
$\Pcal_\Sigma$ is ordered by extension.
I will now prove that $\Pcal_\Sigma$ is proper.
Notice that if this is the case then the sets
$\Dcal_{\alpha} = \{p \in \Pcal_\Sigma:\alpha \in \dom(p)\}$
must be dense.
This is because
$\Dcal^*_x = \{p \in \Pcal_\Sigma:\exists \nu \in \dom(p) ( x \in p(\nu) )\}$
is clearly dense for all $x$ in $X$ and therefore, after forcing with
$\Pcal_\Sigma$, there is always a surjection from
$\{\alpha:\exists p \in G (\alpha \in \dom(p))\}$
onto the uncountable set $X$.

To see that $\Pcal_\Sigma$ is proper, let $p$ be in $\Pcal_\Sigma$ and
$M$ be an elementary submodel of $H(\lambda)$ for
$\lambda$ sufficiently large such that $\Sigma$, $\Pcal_\Sigma$, $p$, and
$H(|{\Pcal_\Sigma}|^+)$ are all in $M$.
Let $\{D_i:i < \omega\}$ enumerate the dense subsets of $\Pcal_\Sigma$ which are in $M$.
We will now build a sequence of conditions $p_0 \geq p_1 \geq \ldots$ by
recursion.
Set $p_0 = p$ and let $p_i$ be given.
Let $E_i$ be the collection of all intersections of the form
$N = N^* \cap X$
where $N^*$ is a countable elementary submodel of
$H(|{\Pcal_\Sigma}|^+)$ containing $H(\theta)$, $D_i$, $\Pcal_\Sigma$, and $p_i$.
Then $E_i \subseteq [X]^{\omega}$ is a club in $M$.
Since $\Sigma(M \cap H(\theta))$ is open and $M \cap H(\theta)$-stationary,
there is a $N_i$ in $E_i \cap \Sigma(M \cap H(\theta)) \cap M$ and an
$x_i$ in $[N_i]^{< \omega}$ such that $[x_i,N_i] \subseteq \Sigma(M \cap H(\theta))$.
Extend $p_i$ to $$q_i = p_i \cup \{(\zeta_i + 1,\hull(p_i(\zeta_i) \cup x_i))\}$$
where $\zeta_i$ is the last element of the domain of $p_i$ and
$\hull(p_i(\zeta_i)\cup x_i)$ is the Skolem hull taken in $H(\theta)$.
Notice that $q_i$ is in $N^*_i$ since $N^*_i$ contains $p_i$ and $H(\theta)$.
Now, working in $N_i^*$, find an extension $p_{i+1}$ of $q_i$
which is in $N_i^* \cap D_i$.
The key observation here is that everything in the range
of $p_{i+1}$ which is not in the range of $p_i$ has an intersection with
$X$ which is in the interval
$[x_i,N_i]$ and therefore in $\Sigma(M \cap H(\theta))$ by the virtue $x_i,N_i$ witnessing
that $\Sigma(M \cap H(\theta))$ is a neighborhood of $N_i$.
Define $p_\infty(\xi) = p_i(\xi)$ if $\xi \leq \zeta_i$ and
$p_\infty(\sup_i \zeta_i) = M \cap H(\theta)$.
It is easily checked that $p_\infty$ is well defined
and that $p_\infty$ is a condition which is moreover $(M,\Pcal_\Sigma)$-generic.
\end{proof}
\section{$\MRP$ and the continuum.}

In this section we will see that $\MRP$ can be used to code reals
in a way which is somewhat reminiscent of the $\SRP$ style coding
methods (e.g. $\phi_{\AC}$, $\psi_{\AC}$ of \cite{Pmax}).
Before we begin, we first need to introduce some notation.

Fix a sequence $\<C_\xi :\xi \in \lim(\omega_1)\>$
such that $C_\xi$ is a cofinal subset of $\xi$ of ordertype $\omega$ for
each limit $\xi < \omega_1$.
Let $N,M$ be countable sets of ordinals such that
$N \subseteq M$, $\otp(M)$ is a limit, and $\sup(N) < \sup(M)$.
Define $$w(N,M) = |\sup(N) \cap \pi^{-1}_M [C_\alpha]|$$
where $\alpha$ is the ordertype of $M$.
A trivial but important observation is that $w$ is \emph{left monotonic} in the
sense that $w(N_1,M) \leq w(N_2,M)$ whenever $N_1 \subseteq N_2 \subseteq M$
and $\sup(N_2) < \sup(M)$.
Also, if $N \subseteq M$ are countable sets of ordinals with
$\sup(N) < \sup(M)$ and $\iota$ is an order preserving map from
$M$ into the ordinals then
$w(N,M) = w(\iota'' N,\iota''M)$.

We will now consider the following statement about a given subset $A$ of $\omega_1$:
\begin{description}

\item[$\upsilon_{\AC}(A)$]
There is an uncountable $\delta < \omega_2$ and
an increasing sequence $\<N_\xi : \xi < \omega_1\>$ which is
club in $[\delta]^{\omega}$ such
that for all limit $\nu < \omega_1$ there is a $\nu_0 < \nu$ such that if
$\xi$ is in $(\nu_0,\nu)$
then $N_\nu \cap \omega_1 \in A$ is equivalent to
$w(N_\xi \cap \omega_1,N_\nu \cap \omega_1) < w(N_\xi,N_\nu)$.

\end{description}

The formula $\upsilon_{\AC}$ is the assertion that $\upsilon_{\AC}(A)$ holds for all $A \subseteq \omega_1$.

\begin{prop}
$\upsilon_{\AC}$ implies that in $L(\Pcal(\omega_1))$ there is a well ordering
of $\Pcal(\omega_1)/\NS$ in length $\omega_2$.
In particular $\upsilon_{\AC}$ implies both $2^{\omega_1} = \omega_2$ and that
$L(\Pcal(\omega_1))$ satisfies the Axiom of Choice.
\end{prop}

\begin{proof}
For each $[A]$ in $\Pcal(\omega_1)/\NS$ let $\delta_{[A]}$ be the least uncountable
$\delta < \omega_2$ such that there is a club $N_\xi$ $(\xi < \omega_1)$ 
in $[\delta]^{\omega}$ such that
for all limit $\nu < \omega_1$ there is a $\nu_0 < \nu$ such that for all
$\nu_0 < \xi < \nu$ we have that
$$N_\nu \cap \omega_1 \in A \textrm{ iff }
w(N_\xi \cap \omega_1,N_\nu \cap \omega_1) < w(N_\xi,N_\nu).$$
It is easily checked that this definition is independent of the choice of representative of
$[A]$.
It is now sufficient to show that if $A$ and $B$ are subsets of $\omega_1$
and $\delta_{[A]} = \delta_{[B]}$ then $[A] = [B]$.
To this end suppose that $\delta < \omega_2$ is
an uncountable ordinal such that for clubs
$N_\xi^A$ $(\xi < \omega_1)$ and $N_\xi^B$ $(\xi < \omega_1)$ in $[\delta]^{\omega}$ we have
for every limit $\nu < \omega_1$ there is a $\nu_0 < \nu$ such that for all $\nu_0 < \xi < \nu$
we have both
$$w(N^A_\xi \cap \omega_1,N_\nu^A \cap \omega_1) < w(N_\xi^A,N_\nu^A) \textrm{ iff }
\xi \in A$$
$$w(N^B_\xi \cap \omega_1,N_\nu^B \cap \omega_1) < w(N_\xi^B,N_\nu^B) \textrm{ iff }
\xi \in B.$$
Now there is a closed unbounded set $C \subseteq \omega_1$ such that if $\xi$ is in $C$ then
$N^A_\xi = N^B_\xi$.
It is easily seen that if $\nu$ is a limit point of $C$ then
$\nu$ is in $A$ iff $\nu$ is in $B$.
\end{proof}

\begin{prop}
$\upsilon_{\AC}$ implies that $2^\omega = 2^{\omega_1}$.
\end{prop}

\begin{proof}
This is virtually identical to the proof that the statement $\theta_{\AC}$ of
Todor\v{c}evi\'c implies $2^\omega = 2^{\omega_1}$ as proved in \cite{BMM_c}.
\end{proof}

The reason for formulating $\upsilon_{\AC}$ is that it is a consequence of $\MRP$.

\begin{thm} \label{MRP_A}
$\MRP$ implies $\upsilon_{\AC}$.
\end{thm}
This is a consequence of the following fact.

\begin{lem}
If $M$ is a countable elementary submodel of $H((2^{\omega_1})^+)$ then the following
sets are open and $M$-stationary:
$$\Sigma_{<}(M) = \{N \in M \cap [\omega_2]^{\omega}:
w(N \cap \omega_1,M \cap \omega_1) < w(N,M \cap \omega_2)\}$$
$$\Sigma_{\geq}(M) = \{N \in M \cap [\omega_2]^{\omega}:
w(N \cap \omega_1,M \cap \omega_1) \geq w(N,M \cap \omega_2)\}.$$
\end{lem}

To see how to prove Theorem \ref{MRP_A} from the lemma, define
$$\Sigma_A(M) =
\left\{
\begin{array}{ll}
\Sigma_{<}(M) &
\textrm{ if }
M \cap \omega_1 \in A \cr
\Sigma_{\geq}(M) &
\textrm{ if }
M \cap \omega_1 \not \in A. \cr
\end{array}
\right.
$$
Now let $N_\xi^*$ $(\xi < \omega_1)$ be a reflecting sequence for $\Sigma_A$.
Let
$$\delta = \bigcup_{\xi < \omega_1} N_\xi \cap \omega_2$$
$$N_\xi = N_\xi^* \cap \omega_2.$$
It is easily verified that $\delta$ is an ordinal, taken
together with $\<N_\xi : \xi < \omega_1\>$, satisfies the
conclusion of $\upsilon_{\AC}(A)$.
We will now return our attention to the proof of the lemma.

\begin{proof}
To see that $\Sigma_{<}(M)$ is $M$-stationary,
let $E \subseteq [\omega_2]^{\omega}$ be a club in $M$.
By the pigeonhole principle, there is a $\gamma < \omega_1$ such that
$$\{\sup(N):N \in E \mand N \cap \omega_1 \subseteq \gamma\}$$ is unbounded in
$\omega_2$.
By elementarity of $M$ there is a $\gamma < M \cap \omega_1$ such that
$$\{\sup(N):N \in E \cap M \mand N \cap \omega_1 \subseteq \gamma\}$$ is unbounded in
$M \cap \omega_2$.
Pick an $N$ in $E \cap M$ such that
$$w(N \cap \omega_2,M \cap \omega_2) =
|\sup(N) \cap \pi^{-1}_{M \cap \omega_2} [C_{\otp(M \cap \omega_2)}]| >
|C_{M \cap \omega_1} \cap \gamma|.$$
Since $N$ is in $M$ and $N$ is countable, $N \subseteq M$ and $\sup(N) < \sup(M)$.
It now follows from the definition of $\Sigma_{<}(M)$ that $N$ is in
$E \cap \Sigma_{<}(M) \cap M$.

To see that $\Sigma_{<}(M)$ is open, let $N$ be in $\Sigma_{<}(M)$.
If $N$ does not have a last element, let $\xi$ be the least element of $N$ greater than
$$\max(\sup(N) \cap \pi^{-1}_{M \cap \omega_2} [C_{\otp(M \cap \omega_2)}]).$$
If $N$ has a greatest element, set $\xi = \max(N)$.
Define $x$ to be the finite set $(N \cap C_{M \cap \omega_1}) \cup \{\xi\}$.
It is easy to see that
$$w(x \cap \omega_1,M \cap \omega_1) = w(N \cap \omega_1,M \cap \omega_1)$$
$$w(x \cap \omega_2,M \cap \omega_2) = w(N \cap \omega_2,M \cap \omega_2)$$
and hence by left monotonicity of $w$ we have that
$[x,N] \subseteq \Sigma_{<}(M)$.

In order to see that $\Sigma_{\geq}(M)$ is $M$-stationary, let $E \subseteq [\omega_2]^{\omega}$ be a
club in $M$ and let $\gamma < \omega_2$ be uncountable such that
$E \cap [\gamma]^{\omega}$ is a club in $[\gamma]^{\omega}$.
By elementarity of $M$, such a $\gamma$ can be found in $M$.
Working in $M$ it is possible to find an $N$ in $E \cap [\gamma]^{\omega}$ such that
$$|N \cap \omega_1 \cap C_{M \cap \omega_1}| \geq |\gamma \cap \pi^{-1}_{M \cap \omega_2}[C_{\otp(M \cap \omega_2)}]|$$
and $\gamma \cap \pi^{-1}_{M \cap \omega_2}[C_{\otp(M \cap \omega_2)}] \subseteq N$.
Then $N$ is in $E \cap \Sigma_{\geq}(M) \cap M$.
The proof that $\Sigma_{\geq}(M)$ is open is similar to the corresponding proof for $\Sigma_{<}(M)$.
\end{proof}

\section{The Bounded Proper Forcing Axiom and the continuum}

In this section we will see that the Bounded Proper Forcing Axiom implies
$\upsilon_{\AC}$.
Before arguing this, I will first give a little context to the result. 
The Bounded Proper Forcing Axiom is equivalent to the assertion that
$$(H(\omega_2),\in) \prec_{\Sigma_1} V^{\Pcal}$$
for every proper forcing $\Pcal$.
That is if $\phi$ is a $\Sigma_1$-formula with a parameter in $H(\omega_2)$
then $\phi$ is true iff it can be forced to be true by some proper forcing.
The original statement of $\BPFA$ is due to Goldstern and Shelah \cite{BPFA}
and is somewhat different,
though equivalent.
The above formulation and its equivalence to the original is due to Bagaria \cite{BFA_gen_abs}.
The consistency strength of $\BPFA$ is much weaker than that of $\PFA$ --- it is exactly
a $\kappa$ reflecting cardinal (such cardinals can exist in $L$)
\cite{BPFA} \cite{local_reflection}.
It should be noted though that many of the consequences of $\PFA$
--- $\MA_{\omega_1}$, the non-existence of S-space and Kurepa trees, the assertion that
all $\omega_1$-dense sets of reals are order isomorphic ---
are actually consequences of $\BPFA$.

Recently there has been a considerable amount of work on bounded forcing axioms and
well orderings of the continuum.
Woodin was the first to give such a proof from the assumption of Bounded Martin's Maximum
and ``there is a measurable cardinal'' \cite{Pmax}.
The question of whether Bounded Martin's Maximum alone sufficed remained an intriguing question. 
In light of Woodin's result, a reasonable approach was to show that $\BMM$ has
considerable large cardinal strength and use this to synthesize the role of the measurable
cardinal.
Asper\'o showed that under $\BMM$ that the dominating number $\dfrak$ is $\omega_2$
\cite{BMM_d_c}.
Soon after, Todor\v{c}evi\'{c} showed that $\BMM$ implied that $\cfrak = \omega_2$
and that, moreover, there is a well ordering of $\Rbb$ which is $\Sigma_1$-definable
from an $\omega_1$-sequence of reals \cite{BMM_c}.
Very recently Schindler showed that $\BMM$ does have considerable consistency strength\cite{BMM_strong}.
It should be noted, however, that it is still unclear whether the measurable cardinal
can be removed from Woodin's argument.

Now we will see that the Bounded Proper Forcing Axiom is already sufficient to
give a definable well ordering of $\Rbb$ from parameters in $H(\omega_2)$.
Notice that for a fixed $A$, $\upsilon_{\AC}(A)$ is a $\Sigma_1$-sentence
which takes the additional parameter $\<C_\xi:\xi \in \lim(\omega_1)\>$
(in order to define $w$).
Further examination of the above proof reveals that for each $A$, there is a proper
forcing which forces $\upsilon_{\AC}(A)$.
Hence the Bounded Proper Forcing Axiom implies $\upsilon_{\AC}$.
Aspero has noted that, unlike statements such as $\psi_{\AC}$ and
$\phi_{\AC}$, the statement $\upsilon_{\AC}$ can be forced over
any model with an inaccessible cardinal.
 
\section{MRP and $\square(\kappa)$}

Recall the following combinatorial principle, defined for $\kappa$ a
regular cardinal greater than $\omega_1$:
\begin{description}

\item[$\square(\kappa)$]
There is a sequence $\<C_\alpha:\alpha < \kappa\>$ such that:
\begin{enumerate}

\item $C_{\alpha+1} = \{\alpha\}$ and $C_\alpha \subseteq \alpha$ is closed
and cofinal if $\alpha$ is a limit ordinal.

\item If $\alpha$ is a limit point of $C_\beta$ then
$C_\alpha = C_\beta \cap \alpha$.

\item There is no club $C \subseteq \kappa$ such that for all
limit points $\alpha$ in $C$ the equality $C_\alpha = C \cap \alpha$
holds.
\end{enumerate}

\end{description}
In this section we will see that $\MRP$ implies that $\square(\kappa)$
fails for all regular $\kappa > \omega_1$.
To this end, let $\<C_\alpha:\alpha < \kappa\>$ be a $\square(\kappa)$-sequence.
The essence of the theorem is contained in the following lemma.

\begin{lem}
If $M$ is a countable elementary submodel of $H(\kappa^+)$ containing
$\<C_\alpha:\alpha < \kappa\>$ then the set $\Sigma(M)$ of all $N \subseteq M \cap \kappa$
such that $\sup(N)$ is not in $C_{\sup(M \cap \kappa)}$ is open and
$M$-stationary.
\end{lem}

Given the lemma, let $N_\nu$ $(\nu < \omega_1)$ be a reflecting sequence
for $\Sigma$ and set $E = \{ \sup(N_\nu \cap \kappa):\nu < \omega_1\}$.
Then $E$ is closed and of order type $\omega_1$.
Let $\beta$ be the supremum of $E$.
Now there must be a limit point $\alpha$ in $E \cap C_\beta$.
Let $\nu$ be such that $\alpha = \sup(N_\nu \cap \kappa)$.
But now there is a $\nu_0 < \nu$ such that
$\sup(N_\xi \cap \kappa)$ is not in $C_\alpha = C_\beta \cap \alpha$ whenever
$\nu_0 < \xi < \nu$.
This means that $\alpha$ is not a limit point of $E$, a contradiction.

Now let us return to the proof of the lemma.
\begin{proof}
First we will check that $\Sigma(M)$ is open.
To see this, let $N$ be in $\Sigma(M)$.
If $N$ has a last element $\gamma$, then $[\{\gamma\},N] \subseteq \Sigma(M)$.
If $N$ does not have a last element, then, since $C_{\sup(M \cap \kappa)}$ is
closed, there is a $\gamma$ in $N$
such that if $\xi < \sup(N)$ is in $C_{\sup(M \cap \kappa)}$ then $\xi < \gamma$.
Again $[\{\gamma\},N] \subseteq \Sigma(M)$.

Now we will verify that $\Sigma(M)$ is $M$-stationary.
To this end, let $E \subseteq [\kappa]^{\omega}$ be a club in $M$.
Let $S$ be the collection of all $\sup(N)$ such that $N$ is in $E$.
Clearly $S$ has cofinally many limit points in $\kappa$.
If $S \cap M$ is contained in $C_{\sup(M \cap \kappa)}$ then we have
that whenever $\alpha < \beta$ are limit points in $S \cap M$,
$$C_\alpha = C_{\sup(M \cap \kappa)} \cap \alpha$$
$$C_\beta = C_{\sup(M \cap \kappa)} \cap \beta$$
and hence $C_\alpha = C_\beta \cap \alpha$.
But, by elementarity of $M$, this means that for all
limit points $\alpha < \beta$ in $S$,
$C_\alpha = C_\beta \cap \alpha$.
This would in turn imply that
the union $C$ of $C_\alpha$ for $\alpha$ a limit point of $S$
is a closed unbounded set such that $C_\alpha = C \cap \alpha$ for all
limit points $\alpha$ of $C$, contradicting the definition of
$\<C_\alpha:\alpha < \kappa\>$.
Hence there is an $N$ in $E$ such that $\sup(N)$ is not in $C_{\sup(M \cap \kappa)}$.
\end{proof}


\begin{thebibliography}{10}

\bibitem{BMM_d_c}
David Asper\'{o}.
\newblock {B}ounded {M}artin's {M}aximum, $\mathfrak{d}$ and $\mathfrak{c}$.
\newblock preprint, 2003.

\bibitem{BFA_gen_abs}
Joan Bagaria.
\newblock Bounded forcing axioms as principles of generic absoluteness.
\newblock {\em Arch. Math. Logic}, 39(6):393--401, 2000.

\bibitem{PFA:Baumgartner}
James Baumgartner.
\newblock Applications of the {P}roper {F}orcing {A}xiom.
\newblock In K.~Kunen and J.~Vaughan, editors, {\em Handbook of Set-Theoretic
  Topology}. North-Holland, 1984.

\bibitem{reals_isomorphic}
James~E. Baumgartner.
\newblock All {$\aleph \sb{1}$}-dense sets of reals can be isomorphic.
\newblock {\em Fund. Math.}, 79(2):101--106, 1973.

\bibitem{set_theory:Bekkali}
M.~Bekkali.
\newblock {\em Topics in set theory}.
\newblock Springer-Verlag, Berlin, 1991.
\newblock {L}ebesgue measurability, large cardinals, forcing axioms,
  $\rho$-functions, Notes on lectures by {S}tevo {T}odor\v cevi\'c.

\bibitem{Yorkshireman}
Keith~J. Devlin.
\newblock The {Y}orkshireman's guide to proper forcing.
\newblock In {\em Surveys in set theory}, volume~87 of {\em London Math. Soc.
  Lecture Note Ser.}, pages 60--115. Cambridge Univ. Press, Cambridge, 1983.

\bibitem{FMS}
Matthew Foreman, Menachem Magidor, and Saharon Shelah.
\newblock {M}artin's {M}aximum, saturated ideals, and nonregular ultrafilters.
  {I}.
\newblock {\em Ann. of Math. (2)}, 127(1):1--47, 1988.

\bibitem{BPFA}
Martin Goldstern and Saharon Shelah.
\newblock The {B}ounded {P}roper {F}orcing {A}xiom.
\newblock {\em J. Symbolic Logic}, 60(1):58--73, 1995.
GoSh:507.
arXiv:math.LO/9501222.

\bibitem{set_theory:Jech}
Thomas Jech.
\newblock {\em Set theory}.
\newblock Perspectives in Mathematical Logic. Springer-Verlag, Berlin, second
  edition, 1997.

\bibitem{set_theory:Kunen}
Kenneth Kunen.
\newblock {\em An introduction to independence proofs}, volume 102 of {\em
  Studies in Logic and the Foundations of Mathematics}.
\newblock North-Holland, 1983.

\bibitem{square_core}
Ernest Schimmerling and Martin Zeman.
\newblock Square in core models.
\newblock {\em Bull. Symbolic Logic}, 7(3):305--314, 2001.

\bibitem{BMM_strong}
Ralf Schindler.
\newblock {B}ounded {M}artin's {M}aximum is stronger than the {B}ounded
  {S}emi-proper {F}orcing {A}xiom.
\newblock preprint 2003.
arXiv:math.LO/0305047.

\bibitem{proper_forcing}
Saharon Shelah.
\newblock {\em Proper and improper forcing}.
\newblock Springer-Verlag, Berlin, second edition, 1998.

\bibitem{iterated_cohen_ext}
Robert Solovay and S.~Tennenbaum.
\newblock Iterated {C}ohen extensions and {S}ouslin's problem.
\newblock {\em Ann. of Math.}, 94:201--245, 1971.

\bibitem{forcing_partition}
Stevo Todor{\v{c}}evi{\'c}.
\newblock Forcing positive partition relations.
\newblock {\em Trans. Amer. Math. Soc.}, 280(2):703--720, 1983.

\bibitem{comparing_continuum}
Stevo Todor{\v{c}}evi{\'c}.
\newblock Comparing the continuum with the first two uncountable cardinals.
\newblock In {\em Logic and scientific methods (Florence, 1995)}, pages
  145--155. Kluwer Acad. Publ., Dordrecht, 1997.

\bibitem{partition_problems}
Stevo Todor\v{c}evi\'{c}.
\newblock {\em Partition Problems In Topology}.
\newblock Amer. Math. Soc., 1989.

\bibitem{local_reflection}
Stevo Todor\v{c}evi\'{c}.
\newblock Localized reflection and fragments of {PFA}.
\newblock In {\em Logic and scientific methods}, volume 259 of {\em DIMACS Ser.
  Discrete Math. Theoret. Comput. Sci.}, pages 145--155. AMS, 1997.

\bibitem{BMM_c}
Stevo Todor\v{c}evi\'{c}.
\newblock Generic absoluteness and the continuum.
\newblock {\em Mathematical Research Letters}, 9:465--472, 2002.

\bibitem{FA_stationary}
Boban Veli{\v{c}}kovi{\'c}.
\newblock Forcing axioms and stationary sets.
\newblock {\em Adv. Math.}, 94(2):256--284, 1992.

\bibitem{Pmax}
W.~Hugh Woodin.
\newblock {\em The {A}xiom of {D}eterminacy, Forcing Axioms, and the
  Nonstationary Ideal}.
\newblock Logic and its Applications. de Gruyter, 1999.

\end{thebibliography}
\end{document}